\documentclass[a4paper,11pt]{article}
\usepackage{amsthm}
\usepackage{extpfeil}

\newlength{\abstractwidth}
\renewcommand{\thanks}[1]{\footnote{#1}} 

\newcommand{\be}{\begin{equation}}
\newcommand{\bea}{\begin{eqnarray}}
\newcommand{\eea}{\end{eqnarray}} \newcommand{\ee}{\end{equation}}
 
 \def\ba{\begin{eqnarray}}
\def\ea{\end{eqnarray}}

\def\C{{\bf C}}

\def\o{\omega}

\def\w{{\bf w}}

\def\z{{\bf z}}

\def\o{\omega}

\def\al{\alpha}
\def\b{\beta}

\def\d{\delta}

\def\o{\omega}

\def\z{\zeta}

\def\ti{\tilde}

\def\C{{\bf C}}

\def\cO{{\cal O}}

\def\cO{{\cal O}}
\def\cI{{\cal I}}

\def\[{{\bf [}}
\def\]{{\bf ]}}

\def\pl{\partial}

\begin{document}
\newtheorem{theorem}{Theorem}
\newtheorem{proposition}{Proposition}
\newtheorem{lemma}{Lemma}
\newtheorem{corollary}{Corollary}
\newtheorem{definition}{Definition}
\newtheorem{conjecture}{Conjecture}
\newtheorem{example}{Example}
\newtheorem{claim}{Claim}
\newtheorem*{thm2}{Theorem 2}

\begin{centering}
 
\textup{\large\bf EXISTENCE OF APPROXIMATE HERMITIAN-EINSTEIN STRUCTURES ON SEMI-STABLE BUNDLES}

\vspace{5 mm}

\textnormal{\large Adam Jacob}

\begin{abstract}
{\small The purpose of this paper is to investigate canonical metrics on a semi-stable vector bundle $E$ over a compact K\"ahler manifold $X$. It is shown that if $E$ is semi-stable, then Donaldson's functional is bounded from below. This implies
that $E$ admits an approximate Hermitian-Einstein structure, generalizing a classic result of Kobayashi for projective manifolds to the K\"ahler case. As an application some basic properties of semi-stable vector bundles over compact K\"ahler manifolds are established, such as the fact that semi-stability is preserved under certain exterior and symmetric products.  }

\end{abstract}

\end{centering}

\begin{normalsize}

\section{Introduction}

The existence of canonical metrics is a fundamental problem in differential geometry. Given a holomorphic vector bundle $E$ over a compact complex Hermitian manifold $X$, a natural metric one could hope to find is a Hermitian-Einstein metric. Specifically this is a metric $H$ on $E$ whose curvature endomorphism $F$ satisfies the following differential equation:
\be
g^{j\bar k}F_{\bar kj}=\mu I\nonumber
\ee
where $\mu$ is a fixed constant and $g_{\bar kj}$ is a Hermitian metric on $T^{1,0}X$. This problem has been solved, first by Narasimhan and Seshadri in the case of curves $\cite{NS}$, then for algebraic surfaces by Donaldson $\cite{Don1}$, and for higher dimensional K\"ahler manifolds by Uhlenbeck and Yau $\cite{UY}$. Using the $C^0$ estimate of Uhlenbeck and Yau, Simpson solved this equation on Higgs bundles and certain non-compact manifolds $\cite{Simp}$. 
Buchdahl extended Donaldson's result to arbitrary complex surfaces in $\cite{Buch}$, and Li and Yau generalized the Donaldson-Uhlenbeck-Yau theorem to any compact complex Hermitian manifold in $\cite{LY}$.
A detailed account of the case of Gauduchon metrics can also be found in the book 
of L\"ubke and Teleman $\cite{TL}$. Most importantly for the purposes of this paper, in $\cite{BaS}$ Bando and Siu were able to extend this theory to metrics on the locally free part of coherent sheaves.

In all cases, the existence of a Hermitian-Einstein metric requires an algebraic notion of stability. We say $E$ is stable (in the sense of Mumford-Takemoto) if for every proper coherent subsheaf ${\cal F}\subset E$,
\be
\frac{deg({\cal F})}{rk({\cal F})}<\frac{deg(E)}{rk(E)}.\nonumber
\ee
With this definition, any irreducible vector bundle $E$ admits a Hermitian-Einstein metric if and only if it is stable. Furthermore the proof of Simpson, and the proof of Siu in $\cite{Siu}$ rely on the fact that a certain functional is bounded from below. This functional, introduced by Donaldson in $\cite{Don1}$, is defined on K\"ahler manifolds and compares two metrics $H_0$ and $H$ on $E$. We denote it by $M(H_0,H,\o)$, where $\o$ is the K\"ahler form on $X$. For fixed $H_0$, its gradient flow is given by:
\be
\label{Donaldson Flow}
H^{-1}\pl_t H=-(g^{j\bar k} F_{\bar kj}-\mu I),
\ee
and one can see that at a critical point of this functional the metric will be Hermitian-Einstein. With this functional in mind, we can now state the following version of the Donaldson-Uhlenbeck-Yau Theorem:

\begin{theorem}
Let $E$ be an indecomposable vector bundle over a compact K\"ahler manifold $X$. Then the following conditions are equivalent:

\medskip
i) E is stable.

\medskip
ii) For any fixed metric $H_0$ on $E$, the Donaldson functional $M(H_0,H,\o)$ is bounded from below and proper.

\medskip
iii) $E$ admits a Hermitian-Einstein metric.
\end{theorem}
While this is an extremely powerful theorem, it leaves many questions unanswered, in particular: What if the bundle is not stable? Does there exist some sort of canonical metric in this case?

In this direction there has been relatively few results,
although recently this problem has been addressed over surfaces by Daskalopoulos and Wentworth in $\cite{DW}$. They proved that for a vector bundle $E$ over a compact K\"ahler surface, the Yang-Mills flow converges to a new metric off the singular set of the graded Harder-Narasimhan-Seshadri filtration. We denote this filtration as $Gr(E)$, and it has finite singular set since it is a torsion free sheaf over a surface. The main analytic tool used in proving this result is a compactness theorem of Uhlenbeck, which states that a sequence of connections along the Yang-Mills flow converges (after going to a subsequence and away from a singular set) to a Yang-Mills connection on a bundle with possibly a different topology than $E$. Daskalopoulos and Wentworth were able to identify this new bundle as $Gr(E)^{**}$ and the limiting connection as coming from a Hermitian-Einstein metric on each of the stable quotients forming $Gr(E)^{**}$. Thus they were able to verify a conjecture of Bando and Siu in the surface case, that the Yang Mills flow will ``break-up" the connection on an un-stable bundle into Hermitian-Eintein connections on the stable quotients of the Harder-Narasimhan-Seshadri filtration $\cite{BaS}$.

With this result in mind, it would be nice to explore the Bando-Siu conjecture in higher dimensional cases.  In this paper we provide some progress towards the semi-stable case. We say that $E$ is semi-stable if for every proper coherent subsheaf ${\cal F}\subset E$,
\be
\frac{deg({\cal F})}{rk({\cal F})}\leq\frac{deg(E)}{rk(E)}.\nonumber
\ee
Our main result is to show the condition of semi-stability is equivalent to the existence of an approximate Hermitian-Einstein structure, which means for all $\epsilon>0$, there exists a metric $H$ on $E$ with curvature $F$ such that:
\be
\sup_X|g^{j\bar k} F_{\bar kj}-\mu I|_{C^0}<\epsilon.\nonumber
\ee
We state our complete result here:

\begin{theorem}
\label{main theorem}
Let $E$ be a holomorphic vector bundle over a compact K\"ahler manifold $X$. Then the following conditions are equivalent:

\medskip
i) E is semi-stable.

\medskip
ii) For any fixed metric $H_0$ on $E$, the Donaldson functional $M(H_0,H,\o)$ is bounded from below.

\medskip
iii) $E$ admits an approximate Hermitian-Einstein structure.
\end{theorem}

When $X$ is a projective manifold, this theorem was first proven by Kobayashi in $\cite{Kob}$. There Kobayashi also conjectures that the result should be true for general compact K\"ahler manifolds, the main difficulty being finding a proof of the lower bound of $M(H_0,H,\o)$ from semi-stability without using certain algebraic facts. We present such a proof in this paper, and are thus able to extend Kobayashi's theorem to the K\"ahler case.

A particularly important feature of the proof of Theorem $\ref{main theorem}$ is that the {\it analytic} property of the lower boundedness of a functional is deduced directly from the {\it algebraic} property of semi-stability. This may be useful for the analogous question in the problem of 
constant scalar curvature K\"ahler metrics.
The analogue of the Donaldson functional is in this case the Mabuchi K-energy \cite{Ma86}, and several analogues
of Mumford-Takemoto stability have also been introduced, including 
Chow-Mumford stability, K-stability (Tian \cite{T97},
Donaldson \cite{Don4}), uniform K-stability (Szekelyhidi \cite{Sze}),
slope-stability (Ross-Thomas \cite{RT}), b-stability (Donaldson \cite{Don2}),
as well as infinite-dimensional notions (Donaldson \cite{Don7}, Phong-Sturm \cite{PS06, PS08}
and references therein). Donaldson \cite{Don6} has shown that Chow-Mumford stability
implies the lower boundedness of the K-energy. It would be very instructive if similar implications
can be established directly from the other notions of stability.
The lower boundedness of the K-energy is an important geometric property in itself. It implies the vanishing of the Futaki
invariant, and plays an important role in the K\"ahler-Ricci flow (see e.g. \cite{PS06, PSSW}), which is a more non-linear analogue of the gradient flow of the Donaldson functional. 

Another potentially interesting feature of Theorem \ref{main theorem} is its proof: a fundamental use is made of the regularization of sheaves, building on the works of Buchdahl \cite{Buch}
and Bando-Siu \cite{BaS},
and it is likely that such regularizations could be relevant to the Bando-Siu conjecture.

We briefly describe the proof of Theorem  $\ref{main theorem}$. The proof of the lower bound is a direct generalization of Donaldson's proof that $M(H_0,H,\o)$ is bounded from below in the semi-stable case if $X$ is a curve. That proof relies on the fact that for any semi-stable vector bundle $E$, one  can find a destabilizing subbundle $S$ with quotient bundle $Q$ such that $S$ is  stable and $Q$ is semi-stable. Then the functional $M$ on $E$ breaks up into the corresponding Donaldson functionals on $S$ and $Q$. Since $S$ is stable, that piece is bounded from below. $Q$ is semi-stable, and of strictly less rank than $E$, so by induction we can keep going until we have rank one bundles, which are stable and thus the Donaldson functional is bounded from below.

For us the key difficulty is that the destabilizing objects $S$ and $Q$ may not be vector bundles (as in the case of curves), but only torsion free sheaves. Thus the bulk of the work goes into defining the functional and corresponding terms on a torsion-free subsheaf $S$ with quotient $Q$. We view these sheaves as holomorphic vector bundles off their singular locus, and the main difficulty is that the induced metrics on these sheaves blow up or degenerate as we approach the singular set. The key tool to help us through this difficulty is an explicit regularization procedure which generalizes a procedure of Buchdahl (from $\cite{Buch}$). After a finite number of blowups, denoted $\pi:\ti X\longrightarrow X$, we can pull back and alter these subsheaves to get smooth vector bundles $\ti S$ and $\ti Q$ on $\ti X$. Degenerate metrics on $S$ and $Q$ can now be indentified with smooth metrics on $\ti S$ and $\ti Q$, allowing many of the desired terms to be computed in this smooth setting, including the Donaldson functional. It also helps with the induction step since now we break apart the Donaldson functional on $\ti Q$, which is a smooth vector bundle with smooth metric, so we only have to worry about subsheaves of smooth vector bundles, and not subsheaves of torsion free sheaves. The proof also relies heavily on the work of Bando and Siu $\cite{BaS}$, and since we use a different regularization procedure than they used, we find it useful to go over some of the important estimates in our case.

The paper is organized as follows: In section $\ref{setup}$ we describe our setup, define induced metrics and state some preliminary results about induced metrics. In section $\ref{reg}$ we describe our regularization procedure, and show how to compare induced metrics with new smooth metrics on the regularized spaces. We also show how many of the associated curvature terms compare. In section $\ref{DFR}$ we define the Donaldson functional on sheaves and show it is bounded from below with the assumption that $M$ has a lower bound on stable sheaves. In section $\ref{stablesheaves}$, we go over the proof that $M$ is bounded below on stable sheaves. Finally in section $\ref{TMTA}$, we provide some applications of our main theorem. 

\medskip
\begin{centering}
{\bf Acknowledgements}
\end{centering}
\medskip

First and foremost, the author would like to thank his thesis advisor, D.H. Phong, for suggesting this problem and for all his advice and support. The author would also like to thank Thomas Nyberg, Tristan Collins, and Luis Garcia for many enlightening discussions while writing this paper. Finally, the author is grateful to Valentino Tosatti for many helpful comments and for suggesting Corollary $\ref{lastcor}$. This research was supported in part by the National Science Foundation grant DMS-07-57372, as well as grant DMS-1204155. The results of this paper are part of the author's Ph.D. thesis at Columbia University.

\section{Preliminaries}
\label{setup}

We begin with some basic facts about holomorphic vector bundles. We also define induced metrics on subsheaves and quotient sheaves.

Let $E$ be a holomorphic vector bundle over the K\"ahler manifold $X$. Locally the K\"ahler form is given by:  
\be
 \o=\frac{i}{2\pi}\,g_{\bar kj}\,dz^j\wedge d\bar z^k.\nonumber
 \ee
where $g_{\bar kj}$ is a Hermitian metric on the holomorphic tangent bundle $T^{1,0}X$.
Let $\Lambda$ denote the adjoint of wedging with $\o$. If $\eta$ is a $(p+1,q+1)$ form, then $\Lambda\eta$ is a $(p,q)$ form with local coeficients $g^{j\bar k}\eta_{\bar QP\bar kj}$, where $P$ and $Q$ are multi-indices of length $p$ and $q$. The volume form on $X$ is given by $\frac{\o^n}{n!}$. For simplicity we write $\o^n$ for the volume form and denote $\frac{\o^{n-1}}{n-1!}$ by $\o^{n-1}$. One can check that for a $(1,1)$ form $\zeta$, we have $(\Lambda \zeta)\,\o^n=\zeta\wedge\o^{n-1}$. Assume that $E$ carries a smooth Hermitian metric $H$. On a local trivialization, for any section $\phi^\al\in\Gamma(X,E)$ we define the unitary-Chern connection by:
\be
\nabla_{\bar k}\phi^\al=\pl_{\bar k}\phi^\al\qquad{\rm and}\qquad\nabla_j\phi^\al=\pl_j\phi^\al+H^{\al\bar\beta}\pl_jH_{\bar\b\gamma}\phi^\gamma.\nonumber
\ee
The curvature of this connection is an endomorphism valued two form:
\be
F:=F_{\bar kj}{}^\al{}_\gamma\,dz^j\wedge d\bar z^k\nonumber,
\ee
where $F_{\bar kj}{}^\al{}_\gamma=-\pl_{\bar k}(H^{\al\bar\b}\pl_jH_{\bar\b\gamma})$. We can compute the {\it degree} of $E$ as follows:
\be
\label{degree}
deg(E)=\frac{i}{2\pi}\int_X{\rm Tr}(F)\wedge\o^{n-1},
\ee
and since $X$ is K\"ahler this definition is independent of a choice of metric on $E$. In the future for notational simplicity we drop the factor $\frac{i}{2\pi}$ as it plays no role in the arguments presented. We define the {\it slope} of $E$ to be\be
\mu(E):=\dfrac{deg(E)}{rk({E})}.\nonumber
\ee

Given a torsion free subsheaf $S$ of $E$, we can construct the following short exact sequence:
\be
\label{sequence}
0\longrightarrow S\xrightarrow{{\phantom {X}}{f}{\phantom {X}}} E\xrightarrow{{\phantom {X}}p{\phantom {X}}} Q\longrightarrow 0,
\ee
where we assume that the quotient sheaf $Q$ is torsion free (by saturating $S$ if necessary). We define the singular set of $Q$ to be $Z:=\{x\in X\,|\, Q_x$ is not free$\}$. Then on $X\backslash Z$, we can view $\eqref{sequence}$ as a short exact sequence of holomorphic vector bundles.  Here, a smooth metric $H$ on $E$ induces a metric ${J}$ on $S$ and a metric $K$ on $Q$. For sections $\psi,\phi$ of $S$, we define the metric $J$ as follows:
\be
\langle\phi,\psi\rangle_J=\langle f(\phi),f(\psi)\rangle_H.\nonumber
\ee
In order to define $K$ on $Q$, we note that a choice of a metric $H$ on $E$ gives a splitting of $\eqref{sequence}$:
\be
\label{splitting}
0\longleftarrow{S}\xleftarrow{{\phantom {X}}{\lambda}{\phantom {X}}} {E}\xleftarrow{{\phantom {X}}{p^\dagger}{\phantom {X}}} {Q}\longleftarrow0.
\ee
Here $\lambda$ is the orthogonal projection from $E$ onto $S$ with respect to the metric $H$. For sections $v,w$ of $Q$, we define the metric $K$ as:
\be
\langle v,w\rangle_K=\langle p^{\dagger}(v),p^{\dagger}(w)\rangle_H.\nonumber
\ee
\begin{definition}
{\rm On $X\backslash Z$ both $S$ and $Q$ are holomorphic vector bundles. We define an} induced metric {\rm on either $Q$ or $S$ to be one constructed as above.} 
\end{definition}
Once we have sequence $\eqref{splitting}$, the second fundamental from $\gamma\in\Gamma(X,\Lambda^{0,1}\otimes Hom(Q,S))$ is given by:
\be
\gamma= \bar\pl\circ p^\dagger. \nonumber
\ee 
We know that for any $q\in\Gamma(X\backslash Z,Q)$, $\gamma(q)$ lies in $S$ since $p$ is holomorphic and $p\circ p^\dagger=I$, thus $p\,(\bar\pl\circ p^\dagger(q))=0$. Now, because the maps $f$ and $p$ degenerate on $Z$, any induced metric may degenerate as we approach the singular set, causing curvature terms to blow up. However, these singularities are not too bad, and the following proposition tells us what control we can expect. 
\begin{proposition}
\label{L2Form}
The second fundamental form of an induced metric is in $L^{2}$, and we have:
\be
\int_{X\backslash Z}g^{j\bar k}{\rm Tr}\,(\gamma_j^{\dagger}\gamma_{\bar k})\,\o^{n}\leq C.\nonumber
\ee
\end{proposition}
 
We prove this proposition in section $\ref{reg}$. We now recall how the curvature on $E$ decomposes onto subbundles and quotient bundles (for instance see $\cite{GH}$). Let $F$ be the curvature of $H$, and denote the curvature of the induced metric $J$ by $F^S$ and the curvature of the induced metric $K$ by $F^Q$. We have:
\be
F^S=F|_{S}-\gamma^\dagger\wedge\gamma
\ee
and
\be
\label{quotientcurv}
F^Q=F|_{Q}+\gamma^\dagger\wedge\gamma.
\ee
Combining these two formulas with the fact that $F$ is smooth implies the following result:
\begin{proposition}
\label{L1curve}
The curvature of an induced metric is in $L^1$.
\end{proposition}
With this proposition we see that formula $\eqref{degree}$ is well defined for an induced metric, and we use the formula to compute the degree of $S$ and $Q$.

\section{Regularization of sheaves}
\label{reg}
In this section we give a procedure to regularize the short exact sequence $\eqref{sequence}$.
This procedure generalizes a procedure of Buchdahl from $\cite{Buch}$ to the higher dimensional case. The main difference is that we do not attempt to regularize arbitrary torsion free sheaves over a Hermitian manifold, we only address the specific case where we have a subsheaf of a vector bundle $E$. In fact, one can view this procedure as a way to regularize the map $f$ so its rank does not drop, allowing us to define a new holomorphic subbundle and quotient bundle. We go over a simple example first which illustrates many of the key points.

Consider the ideal sheaf $\cI$ of holomorphic functions vanishing at the origin in $\C^2$. We can write it as the following holomorphic quotient:
\be
0\longrightarrow \cO\xrightarrow{{\phantom {X}}{f}{\phantom {X}}} \cO^2\xrightarrow{{\phantom {X}}p{\phantom {X}}} \cI\longrightarrow 0,\nonumber
\ee
where the maps are given in matrix form by:
\be
f=\left(
\begin{array}{c}z_1\\z_2
 \end{array} 
\right)\qquad p=\left(
\begin{array}{cc}-z_2&z_1
 \end{array} 
\right).\nonumber
\ee
We blowup at the origin $\pi:\ti\C^2\longrightarrow\C^2$, and let $D=\pi^{-1}(0)$. Pulling back the short exact sequence we get:
\be
0\longrightarrow \cO\xrightarrow{{\phantom {X}}{\pi^*f}{\phantom {X}}} \cO^2\xrightarrow{{\phantom {X}}\pi^*p{\phantom {X}}} \pi^*\cI\longrightarrow 0\nonumber
\ee
(here we are implicitly using the fact that $\pi^*\cO_{\C^2}\cong\cO_{\ti\C^2}$). $\ti\C^2$ can be covered by two coordinate patches $U_i:=\{z_i\neq0\}$, $i=1,2$. On $U_1$ we have coordinates $w_1=z_1$ and $w_2=\frac{z_2}{z_1}$, and we can write our maps as:
\be
\pi^*f=\left(
\begin{array}{c}1\\w_2
 \end{array} 
\right)w_1\qquad \pi^*p=\left(
\begin{array}{cc}-w_2&1
 \end{array} 
\right)w_1.\nonumber
\ee
Similarly on $U_2$ we define coordinates $\zeta_1=\frac{z_1}{z_2}$ and $\zeta_2={z_2}$. Now we define the map $\ti f:\cO(D)\rightarrow\cO^2$ by $\frac{1}{w_1}\pi^*f$ on $U_1$ and $\frac{1}{\zeta_2}\pi^*f$ on $U_2$. This map has domain $\cO(D)$ since we need to multiply local sections on $U_1$ by $\frac{z_2}{z_1}$ to get a map on $U_2$ with the same image in $\cO^2$. Thus we get a new short exact sequence:
\be
0\longrightarrow \cO(D)\xrightarrow{{\phantom {X}}{\ti f}{\phantom {X}}} \cO^2\xrightarrow{{\phantom {X}}\ti p{\phantom {X}}} \cO^2/\cO(D)\longrightarrow 0,\nonumber
\ee
which we say is regularized since now the rank of $\ti f$ does not drop anywhere. Since we know what $\ti p$ is on each coordinate patch, we can explicitly compute the transition functions of $\cO^2/\cO(D)$ in this construction. Given a section $(\eta_1,\eta_2)$ of $\cO^2$, then on $U_1$ $\ti p((\eta_1,\eta_2))=-\frac{z_2}{z_1}\eta_1+\eta_2$ and on $U_2$ we have $\ti p((\eta_1,\eta_2))=-\eta_1+\frac{z_1}{z_2}\eta_2$. Thus the transition function from $U_1$ to $U_2$ is multiplication by $\frac{z_1}{z_2}$, so in this case $\cO^2/\cO(D)\cong\cO(-D)$. Now the regularized sequence can be expressed as:
\be
0\longrightarrow \cO(D)\xrightarrow{{\phantom {X}}{\ti f}{\phantom {X}}} \cO^2\xrightarrow{{\phantom {X}}\ti p{\phantom {X}}} \cO(-D)\longrightarrow 0.\nonumber
\ee
With this example in mind, we now turn to the general procedure.

Once again consider the short exact sequence over $X$:
\be
0\longrightarrow S\xrightarrow{{\phantom {X}}{f}{\phantom {X}}} E\xrightarrow{{\phantom {X}}p{\phantom {X}}} Q\longrightarrow 0,\nonumber
\ee
with $E$ locally free and $Q$ torsion free. Suppose $S$ has rank $s$, $E$ has rank $r$, and $Q$ has rank $q$.  In section $\ref{setup}$ we defined the singular set $Z$ of $Q$, and off this set we can view this sequence as a short exact sequence of holomorphic vector bundles. After choosing coordinates, off of $Z$ we view $f$ as a $r\times s$ matrix of holomorphic functions with full rank. As one approaches $Z$ the rank of $f$ may drop, and it is exactly this behavior that we need to regularize before we can carry out the analysis in later sections.

Let $Z_k$ be the subset of $Z$ where $rk(f)\leq k$. For the smallest such $k$, at a point we can choose coordinates so that $f$ can be expressed as
\be
f=\left(
\begin{array}{cc}I_k&0\\0&g
 \end{array} 
\right),\nonumber
\ee
where $g$ vanishes identically on $Z_k$. Blowing up along $Z_k$ by the map $\pi:\ti X\longrightarrow X$, on a given coordinate $U$ patch let $w$ define the exceptional divisor. Then the pullback of $f$ can be decomposed as follows:
\be
\label{decomf}
\pi^*f=\left(
\begin{array}{cc}I_k&0\\0&\ti g
 \end{array} 
\right)\left(
\begin{array}{cc}I_k&0\\0&w^aI_{s-k}
 \end{array} 
\right),
\ee
where $a$ is the largest power of $w$ we can pull out of the $\pi^*g$. Denote the matrix on the left as $\ti f$ and the matrix on right as $t$. We would like to define $\ti S$ as the image of the sheaf $S$ under the map $t$. Working off $\pi^{-1}(Z)$, if $V$ is another open set we know $\pi^*S$ is a holomorphic vector bundle with transition functions $\{\Phi_{UV}\}$ so that for a section $\psi^{\rho}$ of $\pi^*S$, 
\be
\psi^\rho{}|_{U}=\Phi_{UV}{}^\rho{}_\gamma\psi^\gamma{}|_{V}.\nonumber
\ee
With this in mind, the transition functions $\{\ti\Phi_{UV}\}$ of $\ti S$ can be expressed as:
\be
\ti\Phi_{UV}{}^\rho{}_\gamma=\frac{(w_U){}^{a_\gamma}}{(w_V){}^{a_\rho}}\Phi_{UV}{}^\rho{}_\gamma.\nonumber
\ee
Here $a_\gamma$ is equal to $0$ if $\gamma\leq k$ or $a$ if $\gamma>k$, and $w_U$ defines the exceptional divisor on $U$ (and $w_V$ for $V$). Although these transition functions may blow up as we approach $\pi^{-1}(Z)$, they are useful in understanding how the map $t$ twists up $S$. Now the map $\ti f$ defines a new holomorphic inclusion of the sheaf $\ti S$ into the bundle $\pi^*E$, with a new quotient $\ti Q$. Of course, the rank of $\ti f$ may still drop, but one of two things has happened. Either $rk(\ti f)>k$ on $\pi^*(Z_k)$, or for all $x\in Z_k$, if $\it m_x$ is the maximal ideal at the point $x$, then the smallest power $p$ such that $\it m_x^p$ sits inside the ideal generated by the vanishing of $\ti g$ is smaller than that of $g$. In either case we have improved the regularity of $f$. Of course, this procedure is done in local coordinates, and since $X$ is compact we know $Z_k$ is covered by a finite number of open sets. After a finite number of blowups we can conclude that $rk(\ti f)>k$ everwhere. Thus we can next blowup along $Z_{k+1}$ and continue this process until the rank of $\ti f$ does not drop.

Once the rank does not drop anywhere we have that the map $\ti f$ is holomorphic on $\ti X$. It defines a holomorphic subbundle $\ti S$ of $\pi^*E$ with holomorphic quotient $\ti Q$. Summing up, we have proven the following proposition:
\begin{proposition}
\label{regularization}
Over a compact manifold $X$, let $S$ be a torsion free subsheaf of $E$ with torsion free quotient.  There exits a finite number of blowups
\be
\ti X_N\xrightarrow{{\phantom {X}}{\pi_N}{\phantom {X}}}\ti X_{N-1}\xrightarrow{{\phantom {X}}{\pi_{N-1}}{\phantom {X}}}\cdots \xrightarrow{{\phantom {X}}\pi_2{\phantom {X}}}\ti X_1\xrightarrow{{\phantom {X}}\pi_1{\phantom {X}}}X,\nonumber
\ee
and maps $f_k$ over $\ti X_k$ with the the following properties:

\medskip
i) On each $\ti X_k$ around a given point there exits coordinates so that if $w$ defines the exceptional divisor, there exits a diagonal matrix of monomials in $w$ (denoted $t$) so that
\be
\pi^*_{k-1} f_{k-1}= f_k\,t.\nonumber
\ee

\medskip
ii) The rank of $f_N$ is constant on $\ti X_N$, thus it defines a holomorphic subbundle of $\pi_N^*\circ\cdots\circ\pi_1^*E$ with a holomorphic quotient bundle. 
\end{proposition}
We note that this procedure is consistent with another viewpoint found in Uhlenbeck and Yau $\cite{UY}$. In their paper they view a torsion free sheaf locally as a rational map from $X$ to the Grassmanian $Gr(s,r)$ (this is our map $f$). By Hironaka's Theorem we know this map can be regularized after a finite number of blowups. We follow our procedure in order to find coordinates which let us keep track of how that map changes at each step, and in doing so we can work out how the induced metrics on $\ti S$ and $\ti Q$ change during each step.

\subsection{Induced metrics on regularizations}

We now compute how induced metrics change during regularization. First we need a good local description of these metrics. Recall the short exact sequence $\eqref{sequence}$. Since we view $S$ and $Q$ as holomorphic vector bundles off $Z$, we consider local trivializations for these bundles away from the singular set, and in these coordinates the map $f$ is a matrix of holomorphic functions. For any section $\phi^\al\in\Gamma(X,S)$, we write $f(\phi)=f^\gamma{}_\al\phi^\al\in\Gamma(X,E)$. The induced metric $J_{\bar\b\al}$ is defined by
\be
{J}_{\bar\b\al}\phi^\al\overline{\phi^\b}=H_{\bar\rho\nu}f^\nu{}_\al\phi^\al\overline{f^\rho{}_\beta\phi^\beta},\nonumber
\ee
so we have
\be
\label{desO}
{J}_{\bar\b\al}:=H_{\bar\rho\nu}f^\nu{}_\al\overline{f^\rho{}_\beta}
\ee
in our local trivialization for $S$.

The induced metric $K_{\bar \b\al}$ is defined similarly. Let $q^\al\in\Gamma(X,Q)$. If we recall the splitting $\eqref{splitting}$, then in local coordinates the metric $K_{\bar\beta\al}$ is given by 
\be
K_{\bar\b\al}q^\al\overline{q^\b}=H_{\bar\rho\nu}\,{p^\dagger}^\nu{}_\al q^\al\,\overline{{p^\dagger}^\rho{}_\beta q^\beta},\nonumber
\ee
so
\be
\label{desH}
K_{\bar\b\al}:=H_{\bar\rho\nu}\,{p^\dagger}^\nu{}_\al \overline{{p^\dagger}^\rho{}_\beta }.
\ee
In many cases it will be easier to work with the projection $\lambda$ as opposed to $p^\dagger$. Using the fact that $p$ is surjective we write $q=p(V)$ for some $V\in\Gamma(X,E)$. Then $p^\dagger(q)=p^\dagger p(V)=(I-\lambda) V$. Thus the formula
\be
|q|^2_K=|(I-\lambda)V|^2_H\nonumber
\ee
describes the metric $K$ along with $\eqref{desH}$. We note that $V$ is not unique, however given another $V'$ such that $p(V')=q$, then $p(V-V')=0$, and since $\eqref{sequence}$ is exact we know $(I-\lambda)(V-V')=0$. This justifies the alternate definition of $K$. 

\begin{proposition}
\label{change}
Consider a single blowup from the regularization procedure $\pi:\ti X\longrightarrow X$. Let $J$ and $K$ be metrics induced by $f$ and $\ti J$ and $\ti K$ be metrics induced by $\ti f$, where $\ti f$ is defined by $\eqref{decomf}$. Then if $w$ locally defines the exceptional divisor, there exists non-negative integers $a_{
\al}$ so that:
\be
\pi^*{J}_{\bar\b\al}=w^{a_\al}\overline{w^{a_\b}}\ti{J}_{\bar\b\al}\qquad \qquad
\pi^*K_{\bar\b\al}=\frac{1}{w^{a_\al}\overline{w^{a_\b}}}\ti K_{\bar\b\al}.\nonumber
\ee 
\end{proposition}

\begin{proof}
By $\eqref{decomf}$ we know how $\pi^*f$ decomposes, thus from $\eqref{desO}$ we can see that:
\be
\pi^*{J}_{\bar\b\al}=\pi^*H_{\bar\rho\nu}\pi^*f^\rho{}_\al\overline{\pi^*f^\nu{}_\beta}=\pi^*H_{\bar\rho\nu}w^{a_\al}\,\ti f^\rho{}_\al\overline{w^{a_\b}\ti f^\nu{}_\beta}=w^{a_\al}\overline{w^{a_\b}}\ti{J}_{\bar\b\al}.\nonumber
\ee 
This tells us how the local description of ${J}_{\bar\b\al}$ changes during each step of the regularization. How $K_{\bar\beta\al}$ changes is a little more difficult to see. We note that at each point in $\ti X$ the projection $\lambda$ from $\pi^* E$ onto the image of $\pi^*f$ is equal to the projection $\ti\lambda$ onto the image of $\ti f$. This follows because the only difference between the matrices $\pi^*f$ and $\ti f$ is multiplication by the diagonal matrix $t$ (from $\eqref{decomf}$), which only changes the length of each column vector, not the span of the columns. Thus for $V\in\Gamma(\ti X,\pi^*E)$, we have
\be
(I-\lambda)(V)=(I-\ti\lambda)(V).\nonumber
\ee
We need a formula for how $p^\dagger$ changes under regularization. First we note that on $Q$ the map $p\circ p^\dagger$ is the identity, so for $q$ a section of $\pi^*Q$ we have:
\be
\pi^*p\,\pi^*p^\dagger (q)=q.\nonumber
\ee
We now write $\pi^*p=\ti w\ti p$, where $\ti w$ is a diagonal matrix given by monomials of sections defining the exceptional divisor. So $\ti w\ti p\, \pi^*p^\dagger (q)=q$, and because $\ti w$ is invertible it follows:
\be
\label{q}
\ti p \,\pi^*p^\dagger (q)=\ti w^{-1}q.
\ee
Now since the metric $\pi^*H$ on $\pi^*E$ gives a splitting of the following sequence:
\be
0\longrightarrow \ti S\xrightarrow{{\phantom {X}}{\ti f}{\phantom {X}}} \pi^*E\xrightarrow{{\phantom {X}}\ti p{\phantom {X}}} \ti Q\longrightarrow 0,\nonumber
\ee
we have a map $\ti p^\dagger:\ti Q\rightarrow\pi^* E$. Applying this map to each side of $\eqref{q}$ we get:
\bea
\ti p^\dagger \ti w^{-1}q&=&\ti p^\dagger\ti p \,\pi^*p^\dagger (q)\nonumber\\
&=&(I-\ti\lambda)\pi^*p^\dagger (q)\nonumber\\
&=&(I-\lambda)\pi^*p^\dagger (q)\nonumber\\
&=&\pi^*p^\dagger (q),\nonumber
\eea
where the last line follows from the fact that $\pi^* p^\dagger$ is already perpendicular to the image of  $\pi^*f$. Thus we have shown $\pi^*p^\dagger=\ti p^\dagger \ti w^{-1}$, and plugging this into the formula for the metric we have:
\bea
\pi^*K_{\bar \b\al}s^\al\overline{s^\b}&=&H_{\bar\nu\rho}\,{\pi^*p^\dagger{}}^\rho{}_\al\, s^\al\,\overline{{\pi^*p^\dagger}{}^\nu{}_\beta \,s^\beta}\nonumber\\
&=&\frac{1}{w^{a_\al}\overline{w^{a_\b}}}H_{\bar\nu\rho}\,{\ti p^\dagger}{}^\rho{}_\al\, s^\al\,\overline{{\ti p^\dagger}{}^\nu{}_\beta \,s^\beta}\nonumber\\
&=&\frac{1}{w^{a_\al}\overline{w^{a_\b}}}\ti K_{\bar \b\al}s^\al\overline{s^\b}.\nonumber
\eea
This completes the proof of the proposition.
\end{proof}

\subsection{Transformation of curvature terms}
\label{SFFP}

Now that we know how induced metrics change after each step in the regularization procedure, we can compute how the associated curvature terms change. In this section all computations are local, and we restrict ourselves to working with the sheaf $Q$ with induced metric $K$, since all computation involving the subsheaf $S$ are similar. From now on let $F$ denote the curvature of $K$. First we compute how the trace of curvature changes under regularization.

\begin{lemma}
\label{TrFchanges}
For a single blowup in the regularization procedure $\pi:\ti X\longrightarrow X$, let $w$ locally define the exceptional divisor. Then the following decomposition holds in the sense of currents:
\be
\pi^*{\rm Tr}(F)=\sum_\al a_{\al}\pl\bar\pl\,{\rm log}\,|w|^2+{\rm Tr}(\ti F).\nonumber
\ee
\end{lemma}
Along the course of proving the lemma we will also give a formula for $\pi^*F$ in terms of $\ti F$.
\begin{proof}
We work in a local trivialization and apply Proposition $\ref{change}$:
\bea
\pi^*F_{\bar kj}{}^\al{}_\b&=&-\pl_{\bar k}(\pi^*K^{\al\bar\gamma}\pl_j\pi^*K_{\bar \gamma\b})\nonumber\\
&=&-\pl_{\bar k}(\ti K^{\al\bar\gamma}w^{a_\al}\overline{ w^{a_\gamma}}\pl_j(\frac{1}{w^{a_\b}\overline{w^{a_\gamma}}}\ti K_{\bar\gamma\b})).\nonumber
\eea
Now since $\overline w^{a_\gamma}$ is anti-holomorphic, it follows that
\bea
\pi^*F_{\bar kj}{}^\al{}_\b&=&-\pl_{\bar k}(\ti K^{\al\bar\gamma}w^{a_\al}\pl_j(\frac{1}{w^{a_\b}}\ti K_{\bar\gamma\b}))\nonumber\\
&=&-\pl_{\bar k}(w^{a_\al}\pl_j(\frac{1}{w^{a_\b}})\ti K^{\al\bar\gamma}\ti K_{\bar \gamma\b}+\frac{w^{a_\al}}{w^{a_\b}}\ti K^{\al\bar\gamma}\pl_j \ti K_{\bar\gamma\b})\nonumber\\
&=&a_\al\pl_j\pl_{\bar k}{\rm log}|w|^2\delta^{\al}{}_\b-\pl_{\bar k}(\frac{w^{a_\al}}{w^{a_\b}}\ti K^{\al\bar\gamma}\pl_j \ti K_{\bar\gamma\b}).\label{formula F}
\eea
We can use this last line as a formula for the transformation of $F$. Taking the trace proves the lemma. \end{proof}
Because we need to deal with the pullback of K\"ahler forms under the blowup map, we extend the definition of degree to include these degenerate metrics. 
\begin{definition}
\label{degdeg}
{\rm Let $E$ be a vector bundle on $\ti X$, where $\ti X$ is given by a blowup map $\pi:\ti X\longrightarrow X$. Let $F^E$ be the curvature of a given metric $H$ on $E$, and let $\o$ be a K\"ahler metric on $X$. Then the} degree {\rm of $E$ with respect to $\pi^*\o$ is given by:}
\be
deg(E,\pi^*\o)=\int_{\ti X}{\rm Tr}(F^E)\wedge\pi^*\o^{n-1}.\nonumber
\ee
\end{definition}
Even though the metric $\pi^*\o$ is degenerate on the exceptional divisor, since $\pi^*\o$ is closed this definition is independent of the choice of metric on $E$. Once again if $Q$ is a torsion free sheaf and the curvature of $Q$ is $L^1$ on the locally free part of $Q$, then this definition extends from vector bundles to torsion free sheaves.

In the following lemma asserts that the degree of a sheaf is constant under our regularization procedure as long as we compute with respect to the correct degenerate metrics. 
\begin{lemma}
\label{deg}
Let $\ti Q$ be the regularization of the sheaf $Q$, and let $\pi_N$ represent the composition of blowups needed in the regularization. The following formula holds:
\be
deg(Q,\o)=deg(\ti Q,\pi_N^*\o).\nonumber
\ee
\end{lemma}
\begin{proof}
By Proposition $\ref{L1curve}$ we see the degree of $Q$ is given by:
\be
deg(Q,\o)=\int_X{\rm Tr}(F)\wedge\o^{n-1}.\nonumber
\ee
We now pullback this quantity by the blowup map and regularize $Q$. During each step in the procedure we have:
\bea
\int_X{\rm Tr}(F)\wedge\o^{n-1}&=&\int_{\ti X}\pi^*{\rm Tr}(F)\wedge\pi^*\o^{n-1}\nonumber\\
&=&\int_{\ti X}(\sum_\al a_{\al}\pl\bar\pl\,{\rm log}\,|w|^2+{\rm Tr}(\ti F))\wedge\pi^*\o^{n-1}\nonumber\\
&=&\int_{\ti X}{\rm Tr}(\ti F)\wedge\pi^*\o^{n-1},\nonumber
\eea
since $\pi^*\o$ becomes degenerate along the support of $\pl\bar\pl\,{\rm log}\,|w|^2$. We continue the regularization procedure and after a finite number of blowups $\ti F$ will be smooth. The integral stays the same after each step.
\end{proof}

Proposition $\ref{L2Form}$ also follows from Lemma $\ref{TrFchanges}$, and we now present the proof:
\begin{proof}[Proof of proposition $\ref{L2Form}$]
To prove this result we show that after each step in the regularization procedure $||\gamma||^2_{L^2(X)}=||\ti\gamma||^2_{L^2(\ti X)}$, thus after a finite number of blowups $||\ti\gamma||^2_{L^2(\ti X)}$ will be an integral on a smooth vector bundle over a compact manifold and thus bounded. 
From $\eqref{quotientcurv}$ it follows that that:
\be
{\rm Tr}\,(\gamma^\dagger\wedge\gamma)={\rm Tr}\,(F)-{\rm Tr}\,((I-\lambda)\circ F^E).\nonumber
\ee
Let $\pi$ be a single blowup in our regularization procedure. Pulling back by $\pi$ we compute:
\bea
||\gamma||^2_{L^2(X)}&=&\int_{\ti X}\pi^*{\rm Tr}\,(\gamma^\dagger\wedge\gamma)\wedge\pi^*\o^{n-1}\nonumber\\
&=&\int_{\ti X}({\rm Tr}\,(\ti F)+\sum_\al a_{\al}\pl\bar\pl\,{\rm log}\,|w|^2)\wedge\pi^*\o^{n-1}\nonumber\\
&&-\int_{\ti X}{\rm Tr}\,((I-\lambda)\circ F^E)\wedge\pi^*\o^{n-1}\nonumber\\
&=&\int_{\ti X}{\rm Tr}\,(\ti F)\wedge\pi^*\o^{n-1}-{\rm Tr}\,((I-\ti\lambda)\circ F^E)\wedge\pi^*\o^{n-1}\nonumber\\
&=&\int_{\ti X}{\rm Tr}\,(\ti\gamma^\dagger\wedge\ti\gamma)\wedge\pi^*\o^{n-1}.\nonumber
\eea
Here we used the fact that the projection $\ti \lambda$ is equal to the projection $\lambda$, which we saw in the proof of Proposition $\ref{change}$. This completes the proof of the proposition. \end{proof}

\section{The Donaldson functional on regularizations}
\label{DFR}
In this section we extend the definition of the Donaldson functional to include metrics on torsion free subsheaves $S$ and $Q$. This definition only works for induced metrics, and does not extend to arbitrary metrics defined on the locally free parts of $S$ and $Q$. First we go over the definition of the Donaldson functional on the vector bundle $E$.

Fix a reference metric $H_0$ on $E$. For any other metric $H$ define the endomorphism $h=H_0^{-1}H$, and let $Herm^+(E)$ denote the space of positive definite hermitian endomorphisms of $E$. For $t\in[0,1]$, consider any path $h_t\in Herm^+(E)$ with $h_0=I$ and $h_1=h$, and let $F_t$ be the curvature of the metric $H_t:=H_0\,h_t$ along the path. Then the Donaldson functional is given by:
\be
M(H_0,H,\o)=\int_0^1\int_X{\rm Tr}(F_th^{-1}_t\pl_t h_t)\wedge\o^{n-1}\,dt-\dfrac{\mu(E)}{vol(X)}\int_X{\rm log\,det}(h_1)\,\o^n.\nonumber
\ee
One can check that this definition is independent of the choice of path (for instance see $\cite{Siu}$). Given a blowup map $\pi:\ti X\longrightarrow X$, one can also define the Donaldson functional on a vector bundle over $\ti X$ by integrating with respect to the degenerate metric $\pi^*\o$. Since $\pi^*\o$ is closed the functional will still be independent of path. We now define the Donaldson functional on the sheaves $S$ and $Q$ as follows:

\begin{definition}
{\rm For a subsheaf $S$ of $E$, we define the} Donaldson functional {\rm on $S$ to be:}
\be
M_S(H_0,H,\o):=M_{\ti S}(\ti J_0,\ti J,\pi^*\o),\nonumber
\ee
{\rm for any regularization $\ti S$. Similarly we define the} Donaldson functional {\rm on the quotient sheaf $Q$ to be:}
\be
M_Q(H_0,H,\o):=M_{\ti Q}(\ti K_0,\ti K,\pi^*\o),\nonumber
\ee
{\rm for the regularization $\ti Q$ corresponding to $\ti S$.}
\end{definition}
Here $M_{\ti S}(\ti J_0,\ti J,\pi^*\o)$ and $M_{\ti Q}(\ti K_0,\ti K,\pi^*\o)$ are the Donaldson functionals for the vector bundles $\ti S$ and $\ti Q$ defined using the degenerate metric $\pi^*\o$. We note that the domains of the functionals $M_S$ and $M_Q$ are metrics on the vector bundle $E$, thus this definition only applies to induced metrics and does not extend to arbitrary metrics on $S$ and $Q$. In the following proposition we show that this definition makes sense.
\begin{proposition}
$M_S$ and $M_Q$ are well defined functionals  for any pair of metrics on $E$ and are  independent of the choice of regularization.
\end{proposition}
\begin{proof}
Since the regularization procedure is not unique, we show the functional gives the same value independent of the sequence of blowups chosen. Once again, we prove this proposition for the quotient sheaf $Q$, as the argument works the same for $S$.

As we have seen, a choice of metrics $H_0$ and $H$ on $E$ induce metrics $K_0$ and $K$ on $Q$. Furthermore if we regularize $Q$ we get corresponding induced metrics $\ti K_0$ and $\ti K$ on $\ti Q$. Set $\ti k=\ti K_0^{-1}\ti K$ as the endomorphism relating these two metrics, and let $\ti k_t$, $t\in[0,1]$, be any path in $Herm^+(\ti Q)$ connecting the identity to $\ti k$. We have defined the Donaldson functional on $Q$ to be the following integral:
\be
\label{MtiQ}
M_Q(H_0,H,\o)=\int_0^1\int_{\ti X}{\rm Tr}(\ti F_t\ti k_t^{-1}\pl_t \ti k_t)\,\pi^*\o^{n-1}dt-\dfrac{\mu(E)}{vol(X)}\int_{\ti X}{\rm log\,det}(\ti k_1)\pi^*\o^n.
\ee
Note that the path $\ti k_t$ defines a path $k_t:=\dfrac{w^{a_\al}}{w^{a_\gamma}}\ti k^\al{}_\gamma$ which is an endomorphism of the quotient sheaf one step back in the regularization procedure. Similarly the metrics $(K_0)_{\bar \beta\al}:=\dfrac{1}{w^{a_\al}\overline{w^{a_\b}}} (\ti K_0)_{\bar\b\al}$ and $K_t:=K_0k_t$ are defined one step back in the regularization procedure. Let $F_t$ be the curvature of $K_t$. Then we can compute using formula $\eqref {formula F}$ to get the following equation of currents:
\bea
{\rm Tr}(F_tk_t^{-1}\pl_tk_t)\wedge\pi^*\o^{n-1}&=&(F_t)^\al{}_\b (k_t^{-1})^\b{}_\gamma (\pl_t k_t)^\gamma{}_\al\wedge\pi^*\o^{n-1}\nonumber\\
&=&-\bar\pl(\dfrac{w^{a_\al}}{w^{a_\b}}\ti K^{\al\bar\nu}\pl\ti K_{\bar\nu\b})\dfrac{w^{a_\b}}{w^{a_\gamma}} \ti k^{-1}{}^\b{}_\gamma\dfrac{w^{a_\gamma}}{w^{a_\al}}\pl_t{\ti k}^\gamma{}_\al\wedge\pi^*\o^{n-1}\nonumber\\
&=&{\rm Tr}\,(\ti F\ti k^{-1}\pl_t\ti k)\wedge\pi^*\o^{n-1},\nonumber
\eea
where the third equality holds since $w$ is holomorphic. Thus the first integral does not change at any step in the regularization procedure and we get the following equality:
\be
\label{MQ}
\int_0^1\int_{\ti X}{\rm Tr}(\ti F_t\ti k_t^{-1}\pl_t \ti k_t)\,\pi^*\o^{n-1}dt=\int_0^1\int_X{\rm Tr}(F_tk_t^{-1}\pl_tk_t)\,\o^{n-1}dt.
\ee
Here the integral on the right is only in terms of the initial induced metrics $K_0$ and $K$, where the path $k_t$ is such that $k_0=I$ and $k_1=K_0^{-1}K$. Since the integral in $\eqref{MtiQ}$ is independent of path, we conclude that the integral in $\eqref{MQ}$ is independent of regularization and depends only on the choice of metrics $H_0$ and $H$ on $E$. We now do the same for the second integral of line $\eqref{MtiQ}$.

It helps to write the formula for $k_1$ in matrix notation $k_1=t^{-1} \ti k_1\,t$, where $t$ is the matrix defined in $\eqref{decomf}$. Thus it is clear that det$(k_1)$=det$(\ti k_1)$ for each blowup in the regularization procedure, so once again we can write
\be
\int_{\ti X}{\rm log\,det}(\ti k_1)\pi^*\o^n=\int_X{\rm log\,det}(k_1)\o^n,\nonumber
\ee
where the integral on the right only depends on $K_0$ and $K$. Thus our definition of the Donaldson functional on $Q$ only depends on the choice of metrics $H_0$ and $H$ on $E$.
\end{proof}

Now that we have this definition, we state a decomposition result which plays a major role in the proof of our main theorem. First we assume that $S$ and $Q$ are genuine holomorphic vector bundles, which have the same slope as $E$. In $\cite{Don1}$ Donaldson proved:
\be
M(H_0,H,\o)=M_S({J}_0,{J},\o)+M_Q(K_0,K,\o)+||\gamma||^2_{L^2}-||\gamma_0||^2_{L^2},\nonumber
\ee
where $M(H_0,H,\o)$ is the Donaldson functional on $E$, and $M_S({J}_0,{J},\o)$ and $M_Q(K_0,K,\o)$ are the corresponding Donaldson functionals on $S$ and $Q$. In fact, we can see right away that this decomposition extends to induced metrics on sheaves. Since $M(H_0,H,\o)=M(\pi^*H_0,\pi^*H,\pi^*\o)$, we can pull back the functional and look at the decomposition onto the regularized vector bundles $\ti S$ and $\ti Q$. We get the following decompoistion:
\be
M(\pi^*H_0,\pi^*H,\pi^*\o)=M_{\ti S}(\ti{J}_0,\ti{J},\pi^*\o)+M_{\ti Q}(\ti K_0,\ti K,\pi^*\o)+||\ti\gamma||^2_{L^2}-||\ti\gamma_0||^2_{L^2}.\nonumber
\ee
Now since the $L^2$ norm of the second fundamental form is independent of regularization we get the following lemma:  
\begin{lemma}
\label{decomp}
Let $S$ be a torsion free subsheaf of $E$ with torsion free quotient $Q$. If $S$, $E$, and $Q$ all have the same slope then we have the following decomposition:
\be
M(H_0,H,\o)=M_S(H_0,H,\o)+M_Q(H_0,H,\o)+||\gamma||^2_{L^2}-||\gamma_0||^2_{L^2}.\nonumber
\ee
\end{lemma}

\subsection{A lower bound for the Donaldson functional}
\label{Mainth}
In this section we prove a lower bound for the Donaldson functional on $E$ under the assumption that $M_S(H_0,H,\o)$ is bounded from below for $S$ stable, a fact we shall prove in the next section. We first define a notion of slope and stability with respect to a degenerate metric, using Definition $\ref{degdeg}$:
\begin{definition}
{\rm Let $B$ be a vector bundle on $\ti X$, where $\ti X$ is given by a blowup map $\pi:\ti X\longrightarrow X$. Then the} slope {\rm of $B$ with respect to $\pi^*\o$ is given by:}
\be
\mu(B,\pi^*\o)=\dfrac{deg(B,\pi^*\o)}{rk(B)}.\nonumber
\ee
\end{definition}
\begin{definition}{\rm We say $B$ is} stable {\rm with respect to $\pi^*\o$ if for all proper torsion free subsheaves ${\cal F}\subset B$, we have}
\be
\mu({\cal F},\pi^*\o)<\mu(B,\pi^*\o).\nonumber
\ee
{\rm We say $B$ is} semi-stable {\rm with respect to $\pi^*\o$ if}
\be
\mu({\cal F},\pi^*\o)\leq\mu(B,\pi^*\o).\nonumber
\ee
\end{definition}
\begin{theorem}
\label{lbss}
If $E$ is a semi-stable vector bundle over $X$ compact K\"ahler, than the Donaldson functional is bounded from below on $E$.
\end{theorem}
\begin{proof}$E$ is a semi-stable vector bundle on $X$, so all destabilizing subsheaves have the same slope as $E$. We restrict ourselves to subsheaves which have torsion free quotients. Choose the one with the lowest rank, which we call $S$. Then $S$ is stable since any proper subsheaf of $S$ would be a subsheaf of $E$ and thus would have lower slope (since $S$ was chosen with minimal rank). The torsion free quotient $Q$ has the same slope as $S$ (and $E$), and is semi-stable. 

We now decompose the Donaldson functional into functionals on $S$ and $Q$ using Lemma $\ref{decomp}$. In the next section we show  $M_S(H_0,H,\o)$ is bounded from below since $S$ is stable. $||\gamma||^2_{L^2}$ is a positive term and $||\gamma_0||^2_{L^2}$ is fixed (having only to do with the fixed initial metric $H_0$), so the only remaining term to check is  $M_Q(H_0,H,\o)$. Since $M_Q(H_0,H,\o)=M_{\ti Q}(\ti K_0,\ti K,\pi^*\o)$ for some regularization $\ti Q$, we choose to show the latter term is bounded from below, which is helpful since $\ti K_0$ and $\ti K$ are now smooth metrics on a holomorphic vector bundle $\ti Q$. We need to show $\ti Q$ is semi-stable with respect to $\pi^*\o$, that way we can continue this process of decomposing the functional and use induction on rank. First we prove a few lemmas.

\begin{lemma}If $E$ is semi-stable with respect to $\o$, then $\pi^*E$ is semi-stable with respect to $\pi^*\o$ on $\ti X$.
\end{lemma}
We note this lemma is false if we use the K\"ahler metic $\o_\epsilon=\pi^*\o+\epsilon\sigma$ on $\ti X$ (where $\sigma$ is the pullback of the Fubini-Study metric on the exceptional divisor times a suitable bump function). It only works since $\pi^*\o$ degenerates (see $\cite{Buch}$).
\begin{proof}[Proof of lemma]
Suppose $\pi^*E$ is not semi-stable with respect to $\pi^*\o$. Then it contains a proper subsheaf $\cal F$ of rank $p<r$ such that $\mu({\cal F},\pi^*\o)>\mu(\pi^*E,\pi^*\o)$ (here $r$ is the rank of $E$). Since $\pi$ is an isomorphism off the exceptional divisor, we have that $\mu(\pi_*{\cal F},\o)>\mu(E,\o),$ which would contradict the fact that $E$ is semi-stable if we can show $\pi_*\cal F$ is a proper subsheaf of $E$. Clearly away from $Z$ this is true, and since it is a set of  codimension $\geq2$, off of $Z$ we can view $\pi_*\cal F$ as a rational map from $X$ into the Grassmanian $Gr(p,r)$ (see $\cite{UY}$). We can extend this rational map over $Z$ since $E$ is locally free, thus $\pi_*\cal F$ is a subsheaf of $E$.
\end{proof}

\begin{lemma}
If $\ti Q$ is a torsion free quotient with the same slope as $\pi^*E$, then $\ti Q$ is semi-stable with respect to $\pi^*\o$.
\end{lemma}
\begin{proof}
Suppose $\cal G$ is a subsheaf of $\ti Q$ with $\mu({\cal G},\pi^*\o)>\mu(\ti Q,\pi^*\o)$. Then since we have the exact sequence
\be
0\longrightarrow{\cal G}\longrightarrow\ti Q\longrightarrow{\ti Q/{\cal G}}\longrightarrow 0,\nonumber
\ee
by $\cite{Kob}$ Lemma (7.3) we know $\mu(\ti Q/{\cal G},\pi^*\o)<\mu(\ti Q,\pi^*\o)=\mu(\pi^*E,\pi^*\o)$. We define ${\cal B}:=Ker(\pi^*E\rightarrow\ti Q/{\cal G})$. Then ${\cal B}$ is included in the following exact sequence:
\be
0\longrightarrow{\cal B}\longrightarrow\pi^* E\longrightarrow{\ti Q/{\cal G}}\longrightarrow 0.\nonumber
\ee
Now once again by $\cite{Kob}$ Lemma (7.3) we see $\mu({\cal B},\pi^*\o)>\mu(\pi^*E,\pi^*\o),$ contradicting the semi-stability of $\pi^*E$.
\end{proof}
So $\ti Q$ is semi-stable with respect to $\pi^*\o$, and we continue this process. Recall that the vector bundle $\ti Q$ has smooth metrics $\ti K$ and $\ti K_0$ induced from $H$ and $H_0$ on $E$. Among all subsheaves of $\ti Q$ with the same slope, let $S_1$ be a subsheaf of minimal rank. Then $S_1$ is stable with quotient $Q_1$, and in the next section we show $M_{S_1}(\ti K_0,\ti K,\pi^*\o)$ is bounded from below. Using Lemma $\ref{decomp}$, we can reduce the problem to showing the Donaldson functional is bounded from below on $Q_1$. Blowing up again $\pi_1:\ti X_1\longrightarrow \ti X$ and constructing the regularization $\ti Q_1$ over $\ti X_1$, we see $\ti Q_1$ is semi-stable with respect to $\pi^*_1\circ\pi^*(\o)$ by the previous two Lemmas. Since $\ti Q_1$ has stricly lower rank than $\ti Q$, after a finite number of steps the process will terminate since all rank one sheaves are stable. This proves the lower bound for $M(H_0,H,\o)$.
\end{proof}

\section{A lower bound for stable sheaves}
\label{stablesheaves}
In this section we show that the Donaldson functional $M_S(H_0,H,\o)$ is bounded from below if $S$ is constructed as in the proof of Theorem $\ref{lbss}$. This result relies heavily on $\cite{BaS}$, in which Bando and Siu prove that any stable sheaf admits a Hermitian-Einstein metric off its singular locus. However, the lower bound on the Donaldson functional is not a consequence of this result, but rather it is an essential step in the proof. This is important from the point of view of this paper, since we want the proof of the main theorem to only rely on stability conditions, and not on the existence of any canonical metric. Now, because we use a different regularization procedure than the procedure described in $\cite{BaS}$, we choose to go over parts of the proof here in order to confirm that the necessary details carry over in our case. Furthermore, our proof of the lower bound is different, especially in the induction step used in the proof of Theorem $\ref{lbss}$.

At this point we have only defined the functional $M_S(H_0,H,\o)$ for induced metrics on $S$ (see section $\ref{DFR}$). However, showing this functional is bounded from below is by definition equivalent to showing $M_{\ti S}(\ti J_0,\ti J,\pi^*\o)$ is bounded from below for smooth metrics $\ti J_0$ and $\ti J$. We have that the functional $M_{\ti S}$ is actually defined for any pair of metrics on $\ti S$, allowing us to evolve $\ti J$ by the gradient flow of $M_{\ti S}$. Assume that $\o$ is normalized so that $\int_X\o^n=1$.

As a first step to defining the gradient flow we compute the Euler-Lagrange equation for $M_{\ti S}$. First we only consider a single blowup, and towards the end of the section we consider the case when we have a finite number of blowups. Consider the fixed metric $\ti J_0$ and suppose we have a one parameter family of metrics $\ti J_s$ with $\ti J_1=\ti J$. Since $M_{\ti S}$ is defined via integration along a path and the integral is path independent, assume we are integrating along the path $\ti h_s\in Herm^+(\ti X,\ti S)$ which corresponds to $\ti J_s=\ti J_0\ti h_s$. Let $\ti F_s$ be the curvature of $\ti J_s$. The variation of $M_{\ti S}(\ti J_0,\ti J,\pi^*\o)$ is given by:
\be
\d M_{\ti S}(\ti J_0,\ti J,\pi^*\o)= \int_{\ti X}{\rm Tr}((\Lambda_0\ti F_1-\mu(S,\o) I)\ti h_1^{-1}\pl_s\ti h_1)\pi^*\o^n\nonumber
\ee
(For details of this computation see $\cite{Siu}$). Here $\Lambda_0$ refers to the trace with respect to the degenerate metric $\pi^*\o$. Thus at a critical point of $M$ we have $\Lambda_0 \ti F=\mu(S,\o)I$. We can now consider the flow of metrics on $\ti S$ given by:
\be
\label{degflow}
\ti J^{-1}_t\pl_t\ti J_t=-\Lambda_0\ti F_t+\mu(S,\o) I.
\ee
If $J_t$ is any solution to this flow and we define $M_{\ti S}(t)=M_{\ti S}(\ti J_0,\ti J_t,\pi^*\o)$, then
\be
\pl_tM_{\ti S}(t)=-\int_{\ti X}{\rm Tr}((\Lambda_0\ti F-\mu(S,\o) I)^2)\pi^*\o^n\leq 0,\nonumber
\ee
which is clearly decreasing. Since the flow decreases the value of  $M_{\ti S}$, if it is bounded below along the flow for any initial metric it is bounded from  below in general. A priori it is not clear that the degenerate flow $\eqref{degflow}$ is well defined since $\pi^*g^{j\bar k}$ blows up along the exceptional divisor. Thus our first step is to  prove existence of a solution.

\begin{theorem}[Bando, Siu]
\label{BBAS}
Let $\ti S$ be a vector bundle over $\ti X$, where $\pi:\ti X\rightarrow X$ is the blowup of the K\"ahler manifold $X$ along an analytic subvariety. Let $\pi^*\o$ be the degenerate K\"ahler metric pulled back from $X$. Then there exists a metric $H_0$ and a family of metrics $H(t)$ on $\ti S$ such that $H(0)=H_0$ and $H(t)$ satisfies $\eqref{degflow}$.
\end{theorem}

We prove existence by showing the flow is in fact the limit of existing flows. Let $\pi:\ti X\longrightarrow X$ be the blowup of $X$ on which we construct $\ti S$. On $\ti X$, define the metric $\o_0:=\pi^*\o$ . This metric is degenerate along the exceptional divisor, so we adjust it by adding on a small bump function times the pullback of Fubini-Study metric from the exceptional divisor, which we call $\sigma$ (for details see $\cite{BaS}$,$\cite{Buch}$,$\cite{Buch2}$). This gives us a family of K\"ahler forms $\o_\epsilon$ on $\ti X$, given by $\o_\epsilon=\o_0+\epsilon \sigma$.
Consider $g_\epsilon^{j\bar k}$, which is the inverse of the metric associated to $\o_\epsilon$. Since $\o_\epsilon$ no longer degenerates we know that $g_\epsilon^{j\bar k}$ is smooth. We now can define the standard Donaldson heat flow on $\ti S$ with respect to this new base metric $\o_\epsilon$. We prove uniform bounds in $\epsilon$, showing we can take a subsequence as $\epsilon\longrightarrow 0$ which converges to our degenerate flow $\eqref{degflow}$.

We start out by defining an appropriate conformal change. Set $J_{\epsilon,0}=e^{\phi_\epsilon}\ti J$, where $\phi_\epsilon$ is defined by the equation
 \be
 \Delta_\epsilon\phi_\epsilon={\rm Tr}(-\Lambda_\epsilon\ti F+\mu(\ti S,\o_\epsilon) I).\nonumber
 \ee
This equation admits a smooth solution for $\epsilon >0$ since the right hand side integrates to zero against the volume form $\o_\epsilon^n$. With these initial starting metrics $J_{\epsilon,0}$, the family of flows is given by:
 \be
 \label{thefamily}
 J_{\epsilon,t}^{-1}\pl_t J_{\epsilon,t}=-\Lambda_\epsilon F_{\epsilon,t}+\mu(\ti S,\o_\epsilon) I.
\ee
 
These flows give a family of metrics $J_{\epsilon, t}$ that depend on $\epsilon$ and $t$. As before, we let $h_{\epsilon,t}=J_{\epsilon,0}^{-1}J_{\epsilon,t}$. From this point on we may from time to time drop the subscripts on $J_{\epsilon,t}$ and refer to the metric simply as $J$. To show these flows converge along a subsequence we need the following uniform bounds for the full curvature tensor independent of $\epsilon$:
 \be
 ||F_{\epsilon,t}||_{C^k}\leq C,
 \ee
 for all $k$. This is possible when $0<t_1\leq t\leq t_2<\infty$. In fact, we cannot do better then this, and the bounds fall apart if we send $t_1$ to zero or $t_2$ to infinity. Our first step is an $L^1$ bound.
 \begin{proposition}
 \label{L1Fbound} For all time $t\geq0$, we have the estimate
 \be
  ||\Lambda_\epsilon F_{\epsilon,t}||_{L^1}\leq C,
 \nonumber\ee
 independent of $\epsilon$ and $t$.
 \end{proposition}
 
\begin{proof}
To prove the proposition, we need to work out how $\Lambda_\epsilon F_{\epsilon,t}$ and its norms evolve along the flow. We drop subscripts for simplicity. To compute the evolution of curvature, we use a formula from $\cite{Siu}$:
\be
 \label{evolvF}
 \pl_tF_{\bar ml}=\pl_t(F_{\bar ml}-F^0_{\bar ml})=-\pl_t\nabla_{\bar m}(\nabla_lh\,h^{-1})=-\nabla_{\bar m}\nabla_l(h^{-1}\pl_t h).\nonumber
\ee
Plugging our flow into this equation, we see
\be
\pl_t\Lambda F=g^{l \bar m}\pl_tF_{\bar ml}=-g^{l \bar m}\nabla_{\bar m}\nabla_l(J^{-1}\pl_t J)=g^{l \bar m}\nabla_{\bar m}\nabla_l(\Lambda F)=\overline{\Delta}\Lambda F=\Delta\Lambda F.\nonumber
\ee
The last equality holds because we are taking the Laplacian of the specific endomorphism $\Lambda F$. We now compute how the norm squared evolves:
\bea
\pl_t|\Lambda F|^2&=&\langle\pl_t\Lambda F,\Lambda F\rangle+\langle\Lambda F,\pl_t\Lambda F\rangle.\nonumber\\
&=&\langle\Delta\Lambda F,\Lambda F\rangle+\langle\Lambda F,\overline{\Delta}\Lambda F\rangle.\nonumber
\eea
We note that only in this specific case of the flow \eqref{thefamily} acting on $|\Lambda F|^2$ does the time derivative does not affect the norm $\langle\cdot,\cdot\rangle$ in the first equality, since all the corresponding terms cancel. Next, one sees that:
\be
 \Delta|\Lambda F|^2=\langle\Delta\Lambda F,\Lambda F\rangle+\langle\Lambda F,\overline{\Delta}\Lambda F\rangle+|\nabla\Lambda F|^2+|\overline{\nabla}\Lambda F|^2.\nonumber
\ee
Putting these two lines together we see:
\be
\label{Del2}
\pl_t|\Lambda F|^2=\Delta|\Lambda F|^2-|\nabla\Lambda F|^2-|\overline{\nabla}\Lambda F|^2.
\ee
which implies
\be
\pl_t|\Lambda F|^2\leq\Delta|\Lambda F|^2.\nonumber
\ee
Now, we would like to to prove a similar inequality with $|\Lambda F|$ in place of $|\Lambda F|^2$. The only difficulty is that the laplacian of $|\Lambda F|$ may not be well defined where $|\Lambda F|$ vanishes. To get around this we set $|\Lambda F|_\delta=\sqrt{|\Lambda F|^2+\delta}$ for some small $\delta>0$. Then as $\delta$ goes to zero we have we have $|\Lambda F|_\delta\rightarrow| \Lambda F|$ pointwise. Now all derivatives are defined, allowing us to compute
\be
\pl_t|\Lambda F|_\delta=\frac{\pl_t|\Lambda F|^2}{2|\Lambda F|_\d}\nonumber,
\ee
along with
\be
\Delta|\Lambda F|_\delta=\frac{\Delta|\Lambda F|^2}{2|\Lambda F|_\d}-\frac{|\nabla|\Lambda F|^2|^2}{2(|\Lambda F|^2+\d)^{3/2}}.\nonumber
\ee
Applying equation \eqref{Del2}:
\be
\label{Del3}
(\pl_t-\Delta)|\Lambda F|_\d=-\frac{(|\nabla \Lambda F|^2+|\overline{\nabla}\Lambda F|^2)}{2|\Lambda F|_\d}+\frac{|\nabla|\Lambda F|^2|^2}{2(|\Lambda F|^2+\d)^{3/2}}.
\ee
The right hand side above is shown to be non-positive by the following:
\be
|\nabla|\Lambda F|^2|^2=|\langle\nabla\Lambda F,\Lambda F\rangle+\langle\Lambda F,\overline{\nabla}\Lambda F\rangle|^2\leq(|\nabla \Lambda F|^2+|\overline{\nabla}\Lambda F|^2)(|\Lambda F|^2+\d).\nonumber
\ee
Thus we can integrate \eqref{Del3} and then send $\d$ to zero to get:
\be
\int_{\ti X}\pl_t|\Lambda F|\o^n_\epsilon\leq0.\nonumber
\ee
Now, if the $L^1$ norm of $\Lambda_\epsilon F_{\epsilon,t}$ is finite we can pull the derivative out of the integral:
\be
\pl_t\int_{\ti X}|\Lambda F|\,\o_\epsilon^{n}\leq0.\nonumber
\ee
Since for all $\epsilon$ the $L^1$ norm decreases in time, all we need to do is show that the $L^1$ bound for $\Lambda_\epsilon F_{\epsilon,0}$ is independent of $\epsilon$. To see this we note

\be
\Lambda_\epsilon F_{\epsilon,0}=\Delta_\epsilon\phi_\epsilon I+\Lambda_\epsilon\ti F,\nonumber
\ee
so
\bea
\int_{\ti X}|\Lambda_\epsilon F_{\epsilon,0}|\o^n_\epsilon&\leq&\int_{\ti X}|\Delta_\epsilon\phi_\epsilon|\o^n_\epsilon+\int_{\ti X}|\Lambda_\epsilon\ti F|\o^n_\epsilon\nonumber\\
&=&\int_{\ti X}|{\rm Tr}(-\Lambda_\epsilon\ti F+\mu_\epsilon I)|\o^n_\epsilon+\int_{\ti X}|\Lambda_\epsilon\ti F|\o^n_\epsilon\nonumber\\
&\leq&2\int_{\ti X}|\Lambda_\epsilon\ti F|\o^n_\epsilon+C.\nonumber
\eea
Thus to complete the proof we need to show $||\Lambda_{\epsilon}\ti F||_{L^1}$ is bounded independent of $\epsilon$. Since $\Lambda_{\epsilon}\ti F$ is smooth for $\epsilon>0$, if we can show the bound for $\epsilon=0$ (the degenerate case) we will be done.

First we note that Tr$(\Lambda_0\ti F)=\pi^*$Tr$(\Lambda F)$ since Tr$(\ti F)\wedge\o_0^{n-1}=\pi^*$Tr$(F)\wedge\o_0^{n-1}$. Then since $\pi^*$Tr$(F)$ is in $L^1$ by Proposition $\ref{L1curve}$, we have 
\be
\label{trti}
\int_{\ti X}|{\rm Tr}(\Lambda_0\ti F)|\o_{\epsilon}^n\leq C.
\ee
Furthermore since $\ti J$ is induced from a metric $\pi^*H$ on $\pi^*E$, we have by $\eqref{quotientcurv}$
\be
\Lambda_0\ti F=\pi^*(\Lambda F^E{}|_{S})+\Lambda_0(\gamma^\dagger\wedge\gamma).\nonumber
\ee
Now even though the endomorphism $\Lambda_0(\gamma^\dagger\wedge\gamma)$ is unbounded, we do know it is positive. Thus since $\pi^*\Lambda F^E{}|_{S}$ is the pullback of a smooth endomorphism it follows that the eigenvalues of $\Lambda_0\ti F$ are bounded from below. This fact, along with $\eqref{trti}$, give the desired $L^1$ bound for $\Lambda_0\ti F$. Thus the $L^1$ norm of  $\Lambda_{\epsilon}\ti F$ is independent of $\epsilon$.
\end{proof}
 With this uniform $L^1$ bound, we can now get a $L^\infty$ estimate for $\Lambda_\epsilon F_{\epsilon,t}$.
\begin{proposition}
\label{LFLinfty}
For all $t>0$, the following bound is independent of $\epsilon$:
\be
|\Lambda_\epsilon F_{\epsilon,t}|_{L^\infty}\leq C.\nonumber
\ee
\end{proposition}
\begin{proof}
This bound cannot be extended to $t=0$, since in this case we know that $\Lambda_0F_{0,0}$ is not in $L^\infty$. However, for all times $t>0$ we use a heat kernel estimate. We have seen that this endomorphism evolves by a heat equation $\Lambda_{\epsilon}\pl_t F_{\epsilon,t}=\Delta_{\epsilon}\Lambda_{\epsilon} F_{\epsilon,t}$, thus using the heat kernal $\Phi_{\epsilon,t}$ we have
\be
\Lambda_{\epsilon} F_{\epsilon,t}=\int_{\ti X}\Phi_{\epsilon,t}\Lambda_\epsilon F_{\epsilon,0}\,\o^n_{\epsilon}.\nonumber
\ee
Now Proposition 2 of $\cite{BaS}$ gives a $L^\infty$ estimate for the heat kernel independent of $\epsilon$:
\be
0\leq \Phi_{\epsilon,t}\leq C(t^{-n}+1).\nonumber
\ee
Using this estimate, for any $t>0$ we have
\bea
|\Lambda_\epsilon F_{\epsilon,t}|_{L^\infty}&\leq&\int_{\ti X}|\Phi_{\epsilon,t}\Lambda_\epsilon F_{\epsilon,0}|\,\o^n_{\epsilon}\nonumber\\
&\leq& C(t^{-n}+1)\int_{\ti X}|\Lambda_\epsilon F_{\epsilon,0}|\,\o^n_{\epsilon}\nonumber\\
&\leq& C(t^{-n}+1)\nonumber
\eea
since we have a uniform $L^1$ bound. This proves the proposition.

\end{proof}

Our next step is to prove a uniform bound in $\epsilon$ for Tr$(h_{\epsilon,t})$. Once we get this bound, standard theory for the Donaldson heat flow will give us control of $F$ in $C^k$ for all $k$. 

\begin{proposition}
${\rm Tr}(h_{\epsilon,t})$ is bounded independent of $\epsilon$ for all time $t$ with $0<t_1\leq t\leq t_2<\infty$ .
\end{proposition}
\begin{proof}
Dropping subscripts we have that
\be
\pl_t {\rm Tr}(h)= {\rm Tr}(\pl_t h)=- {\rm Tr}(h(\Lambda F-\mu I)).\nonumber
\ee
Since $t\geq t_1>0$, by the previous proposition $|\Lambda F|\leq C$ for some large positive constant $C$. Then
\be
\pl_t  {\rm Tr}(h)\leq C\,  {\rm Tr}(h).\nonumber
\ee
Set $ {\rm Tr}(h)=f$. We have
\be
\pl_t f- Cf\leq0.\nonumber
\ee
multiplying both sides of the equation by $e^{-Ct}$ we get
\be
e^{-Ct}\pl_t f-Ce^{-Ct}\,f=\pl_t(e^{-Ct}f)\leq0.\nonumber
\ee
Integrating both sides gives
\be
e^{-Ct} {\rm Tr}(h_{\epsilon,t})\leq  {\rm Tr}(h_{\epsilon,0})e^{C\cdot0}=1.\nonumber
\ee
Thus
\be
Tr(h_{\epsilon,t})\leq e^{Ct}\leq e^{Ct_2},\nonumber
\ee
which is independent of $\epsilon$ .
\end{proof}

The conformal change we made in defining $J_{\epsilon,0}$ guarantees that det$\,h=1$ along the flow. Since the trace of $h$ is bounded from above it follows that all the eigenvalues of $h$ are bounded away from zero, and thus $h^{-1}$ is in $L^\infty$. At this point standard theory gives the desired $C^{k}$ bounds of the curvature independent of $\epsilon$. After going to a subsequence $\epsilon_i$, $\epsilon_i\longrightarrow 0$ as $i\rightarrow\infty$, the flows converge to a flow $J_{0,t}$ for $t\in[t_1,t_2]$. This flow is the degenerate flow $\eqref{degflow}$ we hoped to define. It is not unique (it may change if we take a different subsequence or if we choose a different time interval $[t_1',t_2']$), however we can still prove long time existence.

\begin{proposition}[Long time existence]
Given $\ti J$ at time $t_0$, once we choose positive times $t_1$ and $t_2$ to get a degenerate flow for $t\in[t_1,t_2]$, we can extend the flow for all time.
\end{proposition}
\begin{proof}
Recall that we choose subsequence of flows as $\epsilon\rightarrow 0$ to define the degenerate flow.  Now choose a sequence of times $\{t_n\}$ going to infinity. We extend the flow to intervals, $[t_1,t_3],...,[t_1,t_n],...$, always taking subsequences of the defining sequence from the previous step. Since the standard Donaldson heat flow exists for all time, each flow $J_{\epsilon,t}$, $\epsilon>0$, exists for all time, and we can continue this process and get a degenerate flow as $t_n$ goes to infinity. 
\end{proof}
Now that we have the degenerate flow defined for one blowup, we briefly describe the case of multiple blowups. Let $\pi_1,...\pi_k$ be the sequence of blowups given in the regularization procedure for $S$. Assume that $\pi_1$ is the first blowup in the procedure, and thus it is on the singular set with highest codimension. On the final blowup, we have the following K\"ahler form:
\be
\o=\o_0+\epsilon_k\sigma_k+...+\epsilon_1\sigma_1.\nonumber
\ee
If we define the Donaldson heat flow with respect to this K\"ahler form, then letting $\epsilon_1$ go to zero will correspond with the previous work in this section. Thus after going to a subsequence we get a smooth flow for times $t\in[t'_1,\infty)$, $t'_1>t_1$, with respect to the metric
\be
\o=\o_0+\epsilon_k\sigma_k+...+\epsilon_2\sigma_2.\nonumber
\ee
We now repeat this process, which is possible since the $L^1$ bound from Proposition $\ref{L1Fbound}$ is independent of all $\epsilon_i$, including $\epsilon_2$. Thus we get the bounds we need to send $\epsilon_2$ to zero, and along a subsequence get a smooth flow for $t\in[t'_2,\infty)$, where $t'_2>t'_1$. This process continues and after a finite number of steps we have the desired degenerate flow defined for $t\in[t'_k,\infty)$. Choose $t^*\in[t'_k,\infty)$, then Theorem $\ref{BBAS}$ is proved by choosing the initial metric $J_{0,t^*}$.

\begin{proposition}
Let $E$ be a semi-stable vector bundle of rank $r$ over $X$, and let $S$ be a subsheaf of minimal rank among all sheaves with the property $\mu(S,\o)=\mu(E,\o)$. Then the Donaldson functional $M_{\ti S}$ on the regularized vector bundle $\ti S$ over $\ti X$ is bounded from below.
\end{proposition}
\begin{proof}
We have previously shown the degenerate flow is defined for all time with initial metric  $J_{0,t^*}$. We now follow the proof of Simpson from $\cite{Simp}$ to show $M_{\ti S}(t)$ is bounded from below. Suppose we choose $S$ as in Theorem $\ref{lbss}$, so that it has minimal rank among all sheaves with the property  $\mu(S)=\mu(E)$. We actually work along a subsequence of times, which we call $t_i$. Denote $h_i:=h_{0,t_i}$ for simplicity, and let $s_i=$log$(h_i)$. We now use a different form of $M_{\ti S}$, introduced by Donaldson in $\cite{Don3}$. Here, by explicit computation along a specific path, one sees that the functional is given by: 
\be
M_{\ti S}(t_i)=\int_{\ti X}{\rm Tr}(F_{0}\,s_i)\o_0^{n-1}+\int_{\ti X}\sum_{\al\gamma}|\bar\pl s_i^{\gamma}{}_\al|^2\frac{e^{\lambda_\gamma-\lambda_\al}-({\lambda_\gamma-\lambda_\al})-1}{({\lambda_\gamma-\lambda_\al})^2}\o_0^n,\nonumber
\ee
where $\lambda_\al$ are the eigenvalues of $s_i$. Now, because $\o$ is degenerate along the exceptional divisor, we consider the pushforward sheaf $\pi_*\ti S$. Recall $\pi$ is an isomorphism off $Z$, thus $\pi_*\ti S$ is a vector bundle on $X\backslash Z$. Since the set $\pi^{-1}Z$ has measure zero the Donaldson functional can now be expressed as:
\be
M_{\ti S}(t_i)=\int_{X\backslash Z}{\rm Tr}(F_{0}\,s_i)\o^{n-1}+\int_{X\backslash Z}\sum_{\al\gamma}|\bar\pl s_i^{\gamma}{}_\al|^2\frac{e^{\lambda_\gamma-\lambda_\al}-({\lambda_\gamma-\lambda_\al})-1}{({\lambda_\gamma-\lambda_\al})^2}\o^n.\nonumber
\ee
Now we can apply the argument of Simpson. His argument works in this case because the non-compact manifold $X\backslash Z$ satisfies all the assumptions Simpson imposes on the base. Also the key assumption on the bundle, that $\Lambda F_0$ is in $L^\infty$, is satisfied by Proposition $\ref{LFLinfty}$ and the fact that $\pi$ is an isomorphsim off $Z$. We assume by contradiction that there do not exit large constants $C_1$, $C_2$ so that the following estimate holds:
\be
\label{simpest}
||s_i||_{L^1}\leq C_1+C_2M_{\ti S}(t_i).
\ee
Then using the blowup argument of Simpson we can construct a proper torsion free subsheaf $\cal F$ of $\pi_*\ti S$, such that $\mu({\cal F},\o)\geq\mu(\pi_*\ti S,\o)$ and $rk({\cal F})< rk(\ti S)$. Denote $rk({\cal F})=p$. Of course, we assumed that $S$ was stable, not $\pi_*\ti S$, so we do not arrive at a contradiction just yet. However, because $\ti S$ is a subbundle of $\pi^*E$, we have $\pi_*\ti S$ and thus $\cal F$ is a subsheaf of $E$ off $Z$. Once again because $Z$ has codimension two we can view $\cal F$ as locally a rational map into the Grassmanian $Gr(p,r)$ and complete this map over $Z$. So  $\cal F$ is a subsheaf of $E$, and since $E$ is semi-stable we know $\mu({\cal F})=\mu(E)$. However $\cal F$ has rank strictly less than $\ti S$ and thus $S$, contradicting our choice of $S$ as the subsheaf of $E$ with the same slope and minimal rank.

With this contradiction inequality $\eqref{simpest}$ follows, and we can conclude: 
\be
M_{\ti S}(\ti J_0,\ti J,\pi^*\o)\geq-\dfrac{C_1}{C_2}.\nonumber
\ee
By definition $M_S(H_0,H)$ is bounded from below as well. 
\end{proof}
As a final step, we need alter the previous proposition so it can be applied to the induction step in the proof of Theorem $\ref{lbss}$.

\begin{proposition}
Let $E$ be a semi-stable vector bundle of rank $r$ over $\ti X$, where $\ti X$ is given by a blowup $\pi:\ti X\longrightarrow X$. Let $S$ be a subsheaf of minimal rank among all sheaves with the property $\mu(S,\pi^*\o)=\mu(E,\pi^*\o)$. Then the functional $M_{S}(H_0,H,\pi^*\o)$ is bounded from below.
\end{proposition}
\begin{proof}
First we construct the regularization $\ti S$ on the blowup $\pi_1:\ti X_1\longrightarrow\ti X$. As before, we have the degenerate flow defined for all time for some initial metric $J_{0,t^*}$ on $\ti S$. Assume that along a subsequence of times estimate $\eqref{simpest}$ does not hold. We view the sheaf ${\pi}_*{\pi_1}_*\ti S$ as a vector bundle on $X\backslash Z$, and just as in the proof of the previous proposition we use the argument of Simpson from $\cite{Simp}$ to construct a proper torsion free subsheaf $\cal F$ of $\pi_*{\pi_1}_*\ti S$ such that $\mu({\cal F},\o)\geq\mu(\pi_*{\pi_1}_*\ti S,\o)$. From this fact we derive our contradiction.

Since $\ti S$ is a holomorphic subbundle of $\pi_1^*E$, it follows that on $X\backslash Z$, $\cal F$ is a subsheaf of $\pi_* E$. Thus we get a map ${\cal F}\longrightarrow \pi_*E$ defined on all of $X$ given by the composition of restriction to $X\backslash Z$ followed by inclusion. It follows that the natural map $\pi^*\pi_*E\longrightarrow E$ gives us a map:
\be
j:\pi^*{\cal F}\longrightarrow\pi^*\pi_*E\longrightarrow E.\nonumber
\ee
Of course this map may not be injective, however if we quotient out by the kernel of $j$, we can construct a proper subsheaf of $E$:
\be
0\longrightarrow \pi^*{\cal F}/Ker(j)\longrightarrow E.\nonumber
\ee
Because $\pi$ is an isomorphsim off $Z$, we see $j$ is injective off $\pi^{-1}(Z)$, so $Ker(j)$ is a torsion sheaf supported on $\pi^{-1}(Z)$. We will arrive at a contradiction if we can show $\mu(\pi^*{\cal F}/Ker(j),\pi^*\o)=\mu(E,\pi^*\o)$, since $rk(\pi^*{\cal F}/Ker(j))<rk(S)$ and $S$ was chosen to be minimal. Consider the short exact sequence:
\be
\label{kerseq}
0\longrightarrow Ker(j)\longrightarrow \pi^*{\cal F}\longrightarrow\pi^*{\cal F}/Ker(j)\longrightarrow 0.
\ee
$Ker(j)$ is a torsion sheaf, so by Proposition 6.14 from $\cite{Kob}$, the determinant line bundle det$Ker(j)$ admits a non-trivial holomorphic section $\zeta$, which can only vanish along the support of $Ker(j)$. Let $V$ be the vanishing locus of $\zeta$. It follows that:
\be
deg(Ker(j),\pi^*\o)=\int_V\pi^*\o^{n-1}=0,\nonumber
\ee
and the integral on the right is equal to zero since $\pi^*\o$ is degenerate along the exceptional divisor (which contains $V$). Thus by $\eqref{kerseq}$ we have $deg(\pi^*{\cal F}/Ker(j),\pi^*\o)=deg(\pi^*{\cal F},\pi^*\o)$, and since both sheaves have the same rank it follows that:
\be
\mu(\pi^*{\cal F}/Ker(j),\pi^*\o)=\mu(\pi^*{\cal F},\pi^*\o)\geq\mu(\ti S,\pi_1^*\pi^*\o)=\mu(S,\pi^*\o)=\mu(E,\pi^*\o).\nonumber
\ee
$E$ is semi-stable with respect to $\pi^*\o$, so $\mu(\pi^*{\cal F}/Ker(j),\pi^*\o)=\mu(E,\pi^*\o)$, and we have our contradiction. We can now conclude: 
\be
M_{\ti S}(\ti J_0,\ti J,\pi_1^*\pi^*\o)\geq-\dfrac{C_1}{C_2}.\nonumber
\ee
By definition $M_S(H_0,H,\pi^*\o)$ is bounded from below as well. 

\end{proof}

\section{The main theorem and applications}
\label{TMTA}
We begin this section by defining an approximate Hermitian-Einstein structure on a holomorphic vector bundle $E$.
\begin{definition}
{\rm We say $E$ admits an} approximate Hermitian-Einstein structure {\rm if for all $\epsilon>0$, there exists a metric $H$ on $E$ with curvature $F$ such that:}
\be
\sup_X|\Lambda F-\frac{\mu(E)}{Vol(X)} I|_{C^0}<\epsilon.\nonumber
\ee
\end{definition}

With this definition, we now prove Theorem $\ref{main theorem}$ as stated in the introduction.
\begin{proof}This theorem is proven in $\cite{Kob}$ in the case where $X$ is a projective algebraic manifold. The only part of that proof which does not extend to the K\"ahler case is the proof that condition $i)$ implies condition $ii)$. This is exactly what we prove for $X$ K\"ahler in sections $\ref {DFR}$ and $\ref{stablesheaves}$. For a proof that condition $ii)$ implies condition $iii)$ and that condition $iii)$ implies condition $i)$ we direct the reader to $\cite{Kob}$.
\end{proof}
We now state the following applications. In each application $X$ is always assumed to be K\"ahler. The proofs of the first four Corollaries can be found in $\cite{Kob}$, chapter $IV$ section $5$, under the assumption that $E$ admits an approximate Hermitian-Einstein structure. We use Theorem $\ref{main theorem}$ to identify existence of an approximate Hermitian-Einstein structure with $E$ semi-stable. We note that Corollaries 2-4 are not original results, however Theorem $\ref{main theorem}$ provides a natural proof of these statements. For example Corollary 2 is also proven in $\cite{A}$.
\begin{corollary}
If $E$ is semi-stable, so is the symmetric tensor product $S^pE$, and the exterior product $\Lambda^pE$.
\end{corollary}
\begin{corollary}
If $E_1$and $E_2$ are semi-stable vector bundles, so is $E_1\otimes E_2$. \end{corollary}
\begin{corollary}
Let $\hat X$ be a finite unramified covering of $X$ with projection $p:\hat X\longrightarrow X$. If $E$ is a semi-stable vector bundle over $X$, then $p^*E$ is a semi-stable vector bundle over $\hat X$. Also if $\hat E$ is a semi-stable vector bundle over $\hat X$, then $p_*\hat E$ is a semi-stable vector bundle over $E$.
\end{corollary}
\begin{corollary}
Let $E$ be a semi-stable vector bundle of rank $r$ over $X$. Then
\be
\int_X((r-1)c_1(E)^2-2r\,c_2(E))\wedge\o^{n-2}\leq0.\nonumber
\ee
\end{corollary}
In $\cite{Don5}$, Donaldson recalls a result of Atiyah and Bott from $\cite{AtB}$ on the Yang-Mills functional and discusses an analogous result for the Calabi functional.  We restate this result of Atiyah and Bott here.
Let $E$ be a vector bundle over a curve $\Sigma$ and consider a flag $\cal F$ of subbundles:
\be
0=E_0\subset E_1\subset\cdots\subset E_q=E.\nonumber
\ee
Define $\cal F$ to be slope decreasing if $\mu(E_1)>\mu(E_2)>...>\mu(E)$. Let $Q_i=E_i/E_{i-1}$, and define
\be
\Phi({\cal F})^2=\sum_{i=0}^q\mu(Q_i)^2rk(Q_i).\nonumber
\ee
 Then Atiyah and Bott prove that for all $\cal F$ slope drecreasing:
\be
\inf_A||F(A)||^2_{L^2}=\sup_{\cal F}\Phi({\cal F})^2.\nonumber
\ee
In fact, the supremum is attained if $\cal F$ is the Harder-Narasimhan filtration of $E$. Now consider $X$ compact K\"ahler of any dimension and normalize $\o$ to have volume one.  From Theorem $\ref{main theorem}$ we can see directly that:
\begin{corollary}
\label{lastcor}
If $E$ is a semi-stable vector bundle over $X$, then
\be
\label{AB}
\inf_A||\Lambda F(A)||^2_{L^2}=\sup_{\cal F}\Phi({\cal F})^2.
\ee
\end{corollary}
The right hand side is given explicitly by $\mu(E)^2rk(E)$ since the Harder-Narasimhan filtration of $E$ is just $E$. The left hand side can be computed using an approximate Hermitian-Einstein structure. Here on the left we use the Hermitian-Yang-Mills functional, which is equivalent to the Yang-Mills functional up to topological terms. It would be interesting to know if formula $\eqref{AB}$ holds for unstable bundles.
 \end{normalsize}

\vspace{10mm}

\begin{centering}

\textnormal{ Department of Mathematics,
Harvard University,
Cambridge, MA 02138\\
e-mail: ajacob@math.harvard.edu}

\end{centering}

\end{document}